\documentclass [12 pt]{article}
\usepackage{amsfonts}
\title{The Lexicographic First Occurrence of a I-II-III-Pattern}
\author{Torey Burton, Anant P.~Godbole and Brett M.~Kindle\\
Department of Mathematics\\
East Tennessee State University\footnote {E-Mail for Corresponding Author: {\tt godbolea@etsu.edu}}}
\begin{document}
\def\p{\mathbb P}
\def\e{{\mathbb E}}
\def\lr{\left(}
\def\rr{\right)}
\def\var{{\rm Var}}
\def\lc{\left\{}
\def\rc{\right\}}
\def\qed{\vbox{\hrule\hbox{\vrule\kern3pt\vbox{\kern6pt}\kern3pt\vrule}\hrule}}
\def\pnd{{\frac{1}{2}{n\choose k}{n-k\choose k}}}
\def\nk{{n\choose k}}
\def\nkk{{n-k\choose k}}
\def\kr{{k\choose r}}
\def\nkr{{n-k\choose k-r}}
\def\nkj{{n-k\choose k-j}}
\def\hmax{{{k\choose \frac{k^2}{n}} {n-k\choose {k - \frac{k^2}{n}}}}}
\def\kj{{k\choose j}}
\def\pndr{{\frac{1}{2}{n\choose k}{k\choose r}{n-k\choose k-r}}} 
\def\thresh{{\sqrt{\frac{2 \nk}{{k\choose r}{n-k\choose k-r}}}}}
\def\thresha{{\sqrt{\frac{2 \nk}{\nkk}}}}
\def\threshb{{\sqrt{2} e^{\frac{k^2}{2 n}}}}
\newtheorem{thm}{Theorem}
\newtheorem{result}{Result}
\newtheorem{lemma}[thm]{Lemma}
\newtheorem{cor}[thm]{Corollary}
\newtheorem{prop}[thm]{Proposition}
\maketitle
\begin{abstract}
Consider a random permutation $\pi\in{\cal S}_n$.  In this paper, perhaps best classified as a contribution to discrete probability distribution theory, we study the {\it first} occurrence $X=X_n$ of a I-II-III-pattern, where ``first" is interpreted in the lexicographic order induced by the 3-subsets of $[n]=\{1,2,\ldots,n\}$.  Of course if the permutation is I-II-III-avoiding then the first I-II-III-pattern never occurs, and thus $\e(X)=\infty$ for each $n$; to avoid this case, we also study the first occurrence of a I-II-III-pattern given a bijection $f:{\bf Z}^+\rightarrow{\bf Z}^+$.
\end{abstract}
\section{Introduction} Consider a random permutation $\pi\in{\cal S}_n$.  In this short note, perhaps best classified as a contribution to discrete probability distribution theory (AMS Subject Classification 60C05), we study the {\it first} occurrence $X=X_n$ of a I-II-III-pattern, defined as follows:  Order the 3-subsets of $[n]=\{1,2,\ldots,n\}$ in the ``obvious" lexicographic fashion
$$\{1,2,3\}<\{1,2,4\}<\{1,2,5\}<\ldots\{1,2,n\}<\{1,3,4\}<\{1,3,5\}<\ldots$$
$$<\{1,n-1,n\}<\{2,3,4\}<\ldots\{n-2,n-1,n\}.$$  
We say that the first I-II-III-pattern occurs at $\{a,b,c\}$ if $\pi(a)<\pi(b)<\pi(c)$ and if $\pi(d)<\pi(e)<\pi(f)$ does not hold for any $\{d,e,f\}<\{a,b,c\}$.  Of course if the permutation is I-II-III-avoiding, which occurs (\cite{bona}) with probability ${{2n}\choose{n}}/(n+1)!$, then $X=\infty$ and the first I-II-III-pattern never occurs.  Consequently $\e(X)=\infty$ for each $n$; to avoid this case, we also simultaneously present results on the first occurrence of a I-II-III-pattern given a bijection $f:{\bf Z}^+\rightarrow{\bf Z}^+$.  
\section{Results}  In what follows, we use the notation $X=abc$ as short for the event $\{X=\{a,b,c\}\}$, and refer to the case of a bijection on ${\bf Z}^+$ as the $n=\infty$ case.
\begin{prop} For each $n\le\infty$, $$\p(X=12r)=\frac{1}{r-1}-\frac{1}{r}.$$
\end{prop}
\noindent{\bf Proof.} Let $\pi(1),\pi(2),\ldots\pi(r)$ be ordered increasingly as $x_1<x_2<\ldots x_r$.  Then, if $\pi(2)=x_{r-1}$ and $\pi(r)=x_r$, we clearly have $X=12r$.  Conversely if the second largest of $\{\pi(1),\pi(2),\ldots\pi(r)\}$ is not in the $2^{\rm nd}$ spot, then either $\pi(2)=x_r$, or $\pi(2)<x_{r-1}$. In the former case, the only way that we can have $X=12s$ is with $s>r$.  If $\pi(2)<x_{r-1}$, there are two possibilities:  If $\pi(1)<\pi(2)$, then $X=12s$ for some $s<r$, and if $\pi(1)>\pi(2)$ then $X\ne 12s$ for any $s$.  Thus we must have $\pi(2)=x_{r-1}$.  Now for $X$ to equal $12r$, we must necessarily have $\pi(r)=x_r$, or else we would have an earlier occurrence of a I-II-III pattern. It now follows that $X=12r$ if and only if $\pi(2)=x_{r-1}$ and $\pi(r)=x_r$ with the other values arbitrary, so that $\p(X=12r)=\frac{(r-2)!}{r!}=\frac{1}{r-1}-\frac{1}{r}$, as desired.  \hfill\qed
\begin{prop}
If $n=\infty$, then $X=1sr$ for some $2\le s<r$ with probability one.
\end{prop}
\noindent{\bf Proof.}  Obvious.  No matter what value $f(1)$ assumes, there is an $s$ with $f(s)>f(1)$ and an $r>s$ with $f(r)>f(s)$.  Let $s_0$, $r_0$ be the smallest such indices; this yields $X=1s_0r_0$.\hfill\qed

Our ultimate goal is to try to determine the entire probability distribution of $X$; for $n=6$, for example, we can check that the ensemble $\{\p(X=abc):1\le a<b<c\le6\}$ is as follows:
\vfill\eject
 \centerline{Table 1}
 \centerline{\it The First Occurrence of a I-II-III Pattern when $n=6$}
$$\vbox{\halign{
\hfil#\hfil&\qquad
\hfil#\hfil&\qquad
\hfil#\hfil\cr
First I-II-III Pattern& Probability& Cumulative Probability\cr
123& 120/720& 0.1666 \cr
124& 60/720& 0.2500 \cr
125& 36/720& 0.3000 \cr
126& 24/720& 0.3333 \cr
134& 50/720& 0.4028 \cr
135& 28/720& 0.4417 \cr
136& 18/720& 0.4667 \cr
145& 26/720& 0.5028 \cr
146& 16/720& 0.5250 \cr
156& 16/720& 0.5472 \cr
234& 48/720& 0.6139 \cr
235& 22/720& 0.6444 \cr
236& 12/720& 0.6611 \cr
245& 24/720& 0.6944 \cr
246& 12/720& 0.7111 \cr
256& 14/720& 0.7306 \cr
345& 24/720& 0.7639 \cr
346& 10/720& 0.7778 \cr
356& 14/720& 0.7972 \cr
456& 14/720& 0.8167 \cr
$\infty$, i.e. never& 132/720& 1.0000 \cr
}}$$
Recall that the {\it median} of any random variable $X$ is any number $m$ such that $\p(X\le m)\ge 1/2$ and $\p(X\ge m)\ge 1/2$.  Now Propositions 1 and 2 together reveal that for $n=\infty$,
$$\p(X\le 134)\ge\sum_{r=3}^\infty\frac{1}{r-1}-\frac{1}{r}=\frac{1}{2}$$
and
$$\p(X\ge134)=1-\sum_{r=3}^\infty\frac{1}{r-1}-\frac{1}{r}=\frac{1}{2},$$
which shows that $X$ has 134 as its unique median.  For finite $n$, however, the median is larger -- Table 1 reveals, for example, that $m=145$ for $n=6$.
\begin{prop}
For $n=\infty$, 
$$\p(X=1sr)=\frac{1}{(s-1)(r-1)r}=\frac{1}{s-1}\lr\frac{1}{r-1}-\frac{1}{r}\rr.$$
\end{prop}
\noindent{\bf Proof.} It can easily be proved, as in Proposition 1 and keeping in mind that $n=\infty$, that $X=1sr$ iff $\pi(s)=x_{r-1},\pi(r)=x_r$; and $\pi(1)=\max_{1\le j\le s-1}\pi(j)$.  It now follows that
$$\p(X=1sr)=\frac{{{r-2}\choose{s-1}}(s-2)!(r-s-1)!}{r!}=\frac{1}{(s-1)(r-1)r},$$
as claimed.\hfill\qed

Proposition 3 provides us with the entire distribution of $X$ when $n=\infty$; note that 
$$\sum_{s=2}^\infty\sum_{r=s+1}^\infty\frac{1}{s-1}\lr\frac{1}{r-1}-\frac{1}{r}\rr=\sum_{s=2}^\infty\frac{1}{s-1}-\frac{1}{s}=1.$$

The probability of the first I-II-III pattern occurring at positions $12r$ is the same for all $n\le\infty$, noting, though, that for finite $n$, $\sum_s\p(X=12s)=\frac{1}{2}-\frac{1}{n}$. There is, however, a subtle and fundamental difference in general between $\p(X=1rs),r\ge3$, when $n=\infty$ and when $n$ is finite.  We illustrate this fact for $\p(X=13r)$ when $n<\infty$.  Recall from the proof of Proposition 3 that for $X$ to equal $13r$ in the infinite case, we had to have $\pi(3)=x_{r-1}, \pi(r)=x_r$, and $\pi(2)<\pi(1)$.  The above scenario {\it will still}, in the finite case, cause the first I-II-III pattern to occur at positions $13r$, but there is another case to consider. If $n=\pi(2)>\pi(1)$ then it is {\it impossible} for $X$ to equal $12s$ for any $s$; in this case we must have $\pi(3)=x_{r-2}$ and $\pi(r)=x_{r-1}$.  The probability of this second scenario is 
$$\frac{(r-3)!}{n(r-1)!}=\frac{1}{n(r-2)(r-1)}.$$  
Adding, we see that
$$\p(X=13r)=\frac{1}{2}\lr\frac{1}{r-1}-\frac{1}{r}\rr+\frac{1}{n}\lr\frac{1}{r-2}-\frac{1}{r-1}\rr,$$
and, in contrast to the $n=\infty$ case where the net contribution of $\p(X=13r; r\ge 4)$ was 1/6, we have
\begin{eqnarray*}\sum_{r=4}^n\p(X=13r)&=&\frac{1}{2}\sum_{r=4}^n\lr\frac{1}{r-1}-\frac{1}{r}\rr+\frac{1}{n}\sum_{r=4}^n\lr\frac{1}{r-2}-\frac{1}{r-1}\rr\\
&=&\frac{1}{6}-\frac{1}{n(n-1)}.
\end{eqnarray*}
The above example illustrates a general fact:
\begin{thm}
For finite $n$, 
$$\p(X=1sr)=\sum_{k=0}^{s-2}\frac{{{s-2}\choose{k}}{{r-k-2}\choose{s-k-1}}(s-k-2)!(r-s-1)!}{n(n-1)\ldots(n-k+1)(r-k)!}.$$
\end{thm}
\noindent{\bf Proof.}  We may have $k$ of the quantities $\pi(2),\pi(3),\ldots,\pi(s-1)$ being greater than $\pi(1)$, where $k$ ranges from 0 to $s-2$.  In this case, these $\pi$s  must equal, from left to right, $(n,n-1,\ldots,n-k+1)$.  Arguing as before, we must have $\pi(s)=x_{r-1-k}$ and $\pi(r)=x_{r-k}$.  The rest of the proof is elementary.\hfill\qed

Unlike the infinite case, $\sum_s\sum_r\p(X=1sr)\ne1$.  So how much {\it is} $\p(X\ge 234)$, or alternatively, how close to unity is
$$\p(X=1sr)=\sum_{s\ge2}\sum_{r\ge s+1}\sum_{k=0}^{s-2}\frac{{{s-2}\choose{k}}{{r-k-2}\choose{s-k-1}}(s-k-2)!(r-s-1)!}{n(n-1)\ldots(n-k+1)(r-k)!}?$$
We obtain the answer in closed form as follows:  Conditioning on $\pi(1)$, we see that $X\ge234$ iff for $j=n,n-1,\ldots,1$, $\pi(1)=j$, and the integers $n,n-1,\ldots,j+1$ appear from left to right in $\pi$.  Summing the corresponding probabilities ${1}/{n}$, ${1}/{n}$,${1}/{(2!n)}$, ${1}/{(3!n)}$, etc yields
\begin{prop}
$$\p(X\ge234)\sim\frac{e}{n}.$$
\end{prop}
\section{Open Problems}
\begin{itemize}
\item Lexicographic ordering is not our only option; in fact it is somewhat unnatural.  Consider another possibility:  What is 
$$\inf\{k:{\rm there\ is\ a\ I-II-III\ pattern\ in\ } (\pi(1),\ldots,\pi(k))\}?$$
This question is not too hard to answer from Stanley-Wilf theory.  Since 
$$\p((\pi(1),\ldots, \pi(k))\ {\rm is\ I-II-III\ free})=\frac{{{2k}\choose{k}}}{(k+1)!},$$
the probability that $k$ is the first integer for which $(\pi(1),\ldots,\pi(k))$ contains a I-II-III pattern is
$$\frac{{{2k-2}\choose{k-1}}}{k!}-\frac{{{2k}\choose{k}}}{(k+1)!}.$$
A more interesting question is the following:  Conditional on the fact that first I-II-III pattern occurs only after the $k$th ``spot" is revealed, what is the distribution of the first 3-subset, interpreted in the sense of this paper, that causes this to occur?    For example, if we let $n=6$, and are told that the first $k$ for which there is a I-II-II pattern in $(\pi(1),\ldots,\pi(k))$ is 5, what is the chance that the first set that causes this to happen is $\{1,2,5\}, \{1,3,5\}, \{1,4,5\}, \{2,3,5\}, \{2,4,5\}$ or $\{3,4,5\}$?
\item Can the results of this paper be readily generalized to other patterns of length 3?  To patterns of length 4?
\item Theorem 4 and Proposition 5 fall short of providing the exact probabilities $\p(X=rst)$ for $r\ge 2$.  Can these admittedly small probabilities be computed exactly or to a high degree of precision?
\item Does the distribution of $X$ consist, as it does for $n=6$, of a series of decreasing segments with the initial probability of segment $j+1$ no smaller than the final probability of segment $j$?
\end{itemize}
 \medskip

\noindent{\bf Acknowledgment}  The research of the second named author was  supported by NSF Grants DMS-0139286 and DMS-0552730.

\end{document}